\documentstyle[12pt]{article}
\title{\bf The Harnack estimate for the Yamabe flow
           on CR manifolds of dimension 3}
\author{ Shu-Cheng Chang and
  Jih-Hsin Cheng }
\date{ }
\begin{document}
\maketitle
\begin{abstract}
We deform the contact form by the amount of the Tanaka-Webster
curvature on a closed spherical $CR$ three-manifold.
We show that if a contact form evolves with free torsion and
positive Tanaka-Webster curvature as initial data, then a certain
Harnack inequality for the Tanaka-Webster curvature holds.  
\end{abstract}



\section{Introduction}
 Let $M$ denote a closed (i.e., compact without boundary) $CR$ manifold. The Yamabe problem is to find a contact form on $M$ with constant
Tanaka-Webster curvature. In serial papers, Jerison and Lee
initiated the study of this problem. ([JL1],[JL2],[JL3]) In
this paper, we study an evolution equation of contact form
so that the solution is expected to converge to a solution of
the Yamabe problem. 

Let $\theta_{(t)}$ denote a family of contact forms on $M$. 
We can associate to it the so called Tanaka-Webster curvature 
([Ta],[We];see also section 2), denoted $W_{(t)}$.
Our evolution equation, the so called (unnormalized) Yamabe flow,
reads as follows:

\begin{equation}
{\partial_t}{\theta_{(t)}}=-2W_{(t)}{\theta_{(t)}}.
\end{equation}

Write ${\theta_{(t)}}=e^{2{\lambda_{(t)}}}{\hat{\theta}}$ with
respect to 
a fixed contact form ${\hat{\theta}}$. Then we can express the
equation (1.1) in $\lambda_{(t)}$:

\begin{equation}
{\partial_t}{\lambda_{(t)}}=-W_{(t)}.
\end{equation}

Since the linearization of $-W$ with respect to $\lambda$ is
a second-order subelliptic operator, the short time solution
and the uniqueness of (1.2) follows from a standard argument.
(we will discuss this and the long time solution elsewhere)
In this paper we will do the Harnack estimate for
$W$. The first step is to obtain a geometric quantity, usually
called the Harnack quantity. The Harnack quantity is a candidate
quantity for us to do the estimate. Let $\nabla_b$, $\Delta_b$,
$<,>_{J,{\theta}}$ denote the subgradient, sublaplacian, and
the Levi form, respectively. (see section 2 for the definitions)
Following the idea of 
Hamilton in [H1], we can "derive" the following Harnack 
quantity: 

\begin{equation}
Z({\theta},{\eta}){\equiv}2{\Delta_b}W+W^{2}+\frac{W}{t}+
<{\nabla_b}W,{\eta}>_{J,{\theta}}+\frac{1}{8}W|{\eta}|^{2}_{J,{\theta}}
\end{equation}
   
\noindent in which $\eta$ is a Legendrian vector field. (see section 3 for
the definition and more details) In section 4, we prove the following
theorem:

\medskip

$\bf{Theorem\:A}$: Let $(M,{\xi},J)$ be a closed spherical
$CR$ 3-manifold. Suppose there is a contact form $\hat{\theta}$
(together with $J$ defining a positive pseudohermitian structure) with
vanishing torsion and positive Tanaka-Webster curvature. Then under
the Yamabe flow (1.1), 

\begin{equation}
Z({\theta},{\eta}){\geq}0 
\end{equation}

\noindent for any Legendrian vector field $\eta$.

\medskip

Integrating (1.4) from time $t_1$ to time $t_2$, we obtain the
following Harnack inequality. (see section 4 for more details)

$\bf{Theorem\:B}$: Suppose we have the same assumptions as in Theorem A.
Then, under the Yamabe flow (1.1), we have, for all points $x_1$, $x_2$
in $M$ and times ${t_1}<{t_2}$, 

\begin{equation}
\frac{W(x_{2},t_{2})}{W(x_{1},t_{1})}{\geq}(\frac{t_{2}}{t_{1}})^{-2}
exp(-\frac{1}{16}L)
\end{equation}

\noindent where 

$$ L=\inf_{\gamma}{\int}_{t_1}^{t_2}|{\dot{\gamma}}|_{J,{\theta_{(t)}}}
^{2}dt
$$

\noindent and the infimum is taken over all Legendrian paths $\gamma$
with ${\gamma}(t_{1})=x_{1}$ and ${\gamma}(t_{2})=x_{2}$.

In section 5, we show by examples that there are contact forms on the standard $CR$
3-sphere with vanishing torsion and nonconstant positive Tanaka-Webster
curvature. We remark that the study here was motivated by the beautiful work of
Richard Hamilton [H2] on the Ricci flow for surfaces and the work of Ben Chow about the Yamabe flow on locally conformally flat manifolds.([Ch])
In [H2], Hamilton proved, among other results, that on the 2-sphere, if
the initial metric has positive curvature, then the solution metric of
the normalized Ricci flow converges to the limiting metric of constant
curvature. One of the important ingredients in Hamilton's proof is the 
Harnack inequality for the evolved curvatures.

\bigskip

\section{\bf{Basics derived from the flow}}
\setcounter{equation}{0}

Let us first review some basic material in $CR$ geometry. (e.g.,[We],[L1]) 
Let $M$ be a closed 3-manifold with an oriented
contact structure $\xi$. 
There always exists a global contact form $\theta$, obtained by patching together local ones with
a partition of unity. The characteristic vector field of $\theta$ is the 
unique vector field $T$ such that ${\theta}(T)=1$ and ${\cal L}_{T}{\theta}
=0$ or $d{\theta}(T,{\cdot})=0$. A $CR$-structure compatible with $\xi$
is a smooth endomorphism $J:{\xi}{\rightarrow}{\xi}$ such that $J^{2}=-
identity$. 
A pseudohermitian structure compatible with $\xi$ is a
$CR$-structure $J$ compatible with $\xi$ together
with a global contact form $\theta$.   

Given a pseudohermitian structure $(J,{\theta})$, 
we can choose a complex vector field $Z_1$, an eigenvector
of $J$ with eigenvalue $i$, and a complex 1-form ${\theta}^1$ such that
$\{ {\theta},{\theta^1},{\theta^{\bar 1}} \}$ is dual to $\{ T,Z_{1},Z_
{\bar 1} \}$. (${\theta^{\bar 1}}={\bar{({\theta^1})}}$,$Z_{\bar 1}=
{\bar{({Z_1})}}$) It follows that $d{\theta}=ih_{1{\bar 1}}{\theta^1}
{\wedge}{\theta^{\bar 1}}$ for some nonzero real function $h_{1{\bar 1}}$.
If $h_{1{\bar 1}}$ is positive,
we call such a pseudohermitian structure
$(J,{\theta})$ positive, and we can choose a $Z_1$
(hence $\theta^1$) such that $h_{1{\bar 1}}
=1$. That is to say

\begin{equation}
d{\theta}=i{\theta^1}{\wedge}{\theta^{\bar 1}}.
\end{equation}

We'll always assume our pseudohermitian structure $(J,{\theta})$ is positive and $h_{1{\bar 1}}=1$ throughout the paper. 
The pseudohermitian connection of $(J,{\theta})$ is the connection
${\nabla}^{{\psi}.h.}$ on $TM{\otimes}C$ (and extended to tensors) given by

\begin{eqnarray}
{\nabla}^{{\psi}.h.}Z_{1}={\omega_1}^{1}{\otimes}Z_{1},
{\nabla}^{{\psi}.h.}Z_{\bar 1}={\omega_{\bar 1}}^{\bar 1}{\otimes}Z_{\bar 1},
{\nabla}^{{\psi}.h.}T=0 \nonumber
\end{eqnarray}

\noindent in which the 1-form ${\omega_1}^{1}$ is uniquely determined
by the following equation with a normalization condition:

\begin{eqnarray}
&d{\theta^1}={\theta^1}{\wedge}{\omega_1}^{1}+{A^1}_{\bar 1}{\theta}
{\wedge}{\theta^{\bar 1}}& \\
& {\omega_1}^{1}+{\omega_{\bar 1}}^{\bar 1}=0. & \nonumber
\end{eqnarray}
 
The coefficient ${A^1}_{\bar 1}$ in (2.2) is called the (pseudohermitian)
torsion. Since $h_{1{\bar 1}}=1$, $A_{{\bar 1}{\bar 1}}=h_{1{\bar 1}}
{A^1}_{\bar 1}={A^1}_{\bar 1}$. And $A_{11}$ is just the complex
conjugate of $A_{{\bar 1}{\bar 1}}$. Differentiating ${\omega_1}^{1}$
gives

\begin{equation}
d{\omega_1}^{1}=W{\theta^1}{\wedge}{\theta^{\bar 1}}
+2iIm(A_{11,{\bar 1}}{\theta^1}{\wedge}{\theta})
\end{equation}

\noindent where $W$ is the Tanaka-Webster curvature. ([We],[Ta])

We can define the covariant differentiations with respect to the
pseudohermitian connection. For instance, $f_{,1}=Z_{1}f$, $f_{1{\bar 1}}
=Z_{\bar 1}Z_{1}f-{\omega_1}^{1}(Z_{\bar 1})Z_{1}f$ for a (smooth)
function $f$.
(see,e.g.,section 4 in [L1]) We define the subgradient operator $\nabla_b$ and the sublaplacian operator $\Delta_b$ by 

\begin{eqnarray}
&{\nabla_b}f=f_{,{\bar 1}}Z_{1}+f_{,1}Z_{\bar 1},&\nonumber\\
&{\Delta_b}f=f_{,1{\bar 1}}+f_{,{\bar 1}1},&\nonumber
\end{eqnarray}

\noindent respectively. (notice the sign difference for $\Delta_b$ 
in [L1]) We also define the Levi form $<,>_{J,{\theta}}$ by

$$
 <V,U>_{J,{\theta}}=2d{\theta}(V,JU)=v_{1}u_{\bar 1}+v_{\bar 1}u_{1}
$$

\noindent for $V=v_{1}Z_{\bar 1}+v_{\bar 1}Z_{1}$,$U=u_{1}Z_{\bar 1}
+u_{\bar 1}Z_{1}$ in $\xi$. (note that the second equality follows
from (2.1) and our definition is different from the one in [L1] by
a factor 2) The associated norm is defined as usual: $|V|_{J,{\theta}}^2
=<V,V>_{J,{\theta}}$. Let ${\hat{\theta}},{\hat{\theta}}^{1},{\hat{\theta}}^{\bar 1}$ satisfy (2.1). Now consider the change of contact form:
${\theta}=e^{2{\lambda}}{\hat{\theta}}$. Choose ${\theta^1}=e^{\lambda}
({\hat{\theta}}^{1}+2i{\lambda_{,{\bar 1}}}{\hat{\theta}})$ such that 
$h_{1{\bar 1}}={\hat h}_{1{\bar 1}}$(=1 by assumption). One checks easily
that ${\theta},{\theta^1},{\theta}^{\bar 1}$ satisfies (2.1). Then the
associated connection form ${\omega_1}^{1}$, torsion $A_{11}$, and
Tanaka-Webster curvature $W$ transform as follows: (cf. section 5 in
[L1])

\begin{eqnarray}
&{\omega_1}^{1}={\hat{({\omega_1}^{1})}}+3({\lambda_{,1}}
{\hat{\theta}}^{1}-{\lambda_{,{\bar 1}}}{\hat{\theta}}^{\bar 1})+
i({\hat{\Delta}}_{b}{\lambda}+4|{\hat{\nabla}}_{b}{\lambda}|_
{J,{\hat{\theta}}}^{2}){\hat{\theta}}&\\
&A_{11}=e^{-2{\lambda}}({\hat A}_{11}+2i{\lambda_{,11}}-4i({\lambda}_
{,1})^{2})& \\
&W=e^{-2{\lambda}}(-4{\hat{\Delta}}_{b}{\lambda}-4|{\hat{\nabla}}_{b}
{\lambda}|^{2}_{J,{\hat{\theta}}}+{\hat{W}}).& 
\end{eqnarray}

\noindent Here the operators or quantities with "hat" are with respect
to the coframe $({\hat{\theta}},{\hat{\theta}}^{1},{\hat{\theta}}^{\bar 1})$, and so are the covariant derivatives of $\lambda$. Now consider
a family of contact forms $\theta_{(t)}=e^{2{\lambda_{(t)}}}
{\hat{\theta}}$, a solution to the Yamabe flow (1.1) or (1.2). 

\medskip

$\bf{Lemma\:2.1.}$ Under the Yamabe flow (1.1)($W=W_{(t)}$ for short), 
we have 

\begin{eqnarray}
&{\dot W}=4{\Delta_b}W+2W^{2}&\\
&{\dot{(A_{11})}}=2WA_{11}-2iW_{,11}&
\end{eqnarray}

\noindent in which $\Delta_b$, the torsion, and covariant derivatives are with respect to $\theta_{(t)}$ and induced coframes as shown previously.

\medskip

$\bf{Proof}$: We will omit the t-dependence for simplicity of
notation if no confusion occurs. First note that $Z_{1}=e^{-{\lambda}}
{\hat Z}_{1}$ and ${\nabla_b}f=e^{-2{\lambda}}{\hat{\nabla}}_{b}f$. It
follows that

\begin{eqnarray}
&{\Delta_b}f=e^{-2{\lambda}}({\hat{\Delta}}_{b}f+2<{\hat{\nabla}}_{b}
{\lambda},{\hat{\nabla}}_{b}f>_{J,{\hat{\theta}}}).&
\end{eqnarray}

Differentiating (2.6) with respect to $t$, we obtain (2.7) by making
use of (1.2) and (2.9). From (2.4), it is easy to see that 

\begin{eqnarray}
&{\hat{({\omega_1}^{1})}}({\hat Z}_{1})=e^{\lambda}({\omega_1}^{1}(Z_{1})-3
{\lambda_{,1}}).&
\end{eqnarray}

Substituting (2.10) in the expression of 
$W_{,{\hat 1}{\hat 1}}$, we obtain

\begin{eqnarray}
&W_{,{\hat 1}{\hat 1}}=e^{2{\lambda}}(W_{,11}+4{\lambda_{,1}}W_{,1}).&
\end{eqnarray}

Now differentiating (2.5) with respect to $t$ and making use of (1.2),(2.11), we finally reach (2.8).

\begin{flushright}
Q.E.D.
\end{flushright}

\medskip

Now applying the maximum principle to (2.7), we obtain 

\medskip

$\bf{Corollary\:2.2}$. Suppose $(M,J,{\theta_{(0)}})$ is closed
with $W{\geq}c>0$. Then the inequality
$W{\geq}c>0$ is preserved under the Yamabe flow (1.1).

\medskip

The following formula will be used to compute the evolution of 
${\Delta_b}W$.

\medskip

$\bf{Lemma\:2.3}$. Under the Yamabe flow (1.1), we have

\begin{eqnarray}
&{\partial_t}({\Delta_b}f)={\Delta_b}({\dot f})+2W{\Delta_b}f
  -2<{\nabla_b}W,{\nabla_b}f>_{J,{\theta}}&
\end{eqnarray}

\noindent for a (smooth) real-valued function $f=f(x,t)$ defined
on $M{\times}R$. (note that we have suppressed the t-dependence in the
above expression)

\medskip

$\bf{Proof}$: Differentiating $Z_{1}=e^{-{\lambda}}{\hat Z}_{1}$
and (2.4) with respect to $t$ gives 

\begin{eqnarray}
&{\dot Z}_{1}=WZ_{1}&\\
&{\dot{({\omega_1}^{1})}}=-3W_{,1}{\theta^1}+3W_{,{\bar 1}}
{\theta^{\bar 1}}\:(mod\:{\theta})&
\end{eqnarray}

\noindent by (1.2). Now our formula (2.12) follows from (2.13),
(2.14) by a straightforward computation.

\begin{flushright}
Q.E.D.
\end{flushright}

\bigskip

\section{\bf{The Harnack quantity}}
\setcounter{equation}{0}

In this section, we apply Hamilton's general method for obtaining
a potential quantity for the Harnack estimate. First we need to
know what the soliton equation is supposed to be for our flow 
(1.1). Let $\phi_t$ be a family of $CR$ automorphisms. Suppose
${\phi}_t^{\star}{\theta_{(t)}}$ converges to a fixed contact
form $\theta$ and differentiating 
${\phi}_t^{\star}{\theta_{(t)}}$ with respect to $t$ converges
to $0$. Then $\theta$ satisfies the following equation:

\begin{eqnarray}
&{\cal L}_{X_f}{\theta}-2W{\theta}=0&
\end{eqnarray}

\noindent in which $X_f$ is a $CR$ vector field parametrized by
a real-valued function $f$. We can write (e.g.,[CL])

\begin{eqnarray}
&X_{f}=-fT+if_{,1}Z_{\bar 1}-if_{,{\bar 1}}Z_{1}.&
\end{eqnarray}

We call (3.1) the soliton equation of the flow (1.1). (A solution
$\theta$ is called a soliton) Substituting (2.1), (3.2) in the 
formula for the Lie derivative, we can reduce (3.1) to
$f_{,0}=-2W$. (recall that $f_{,0}=Tf$ where $T$ is the 
characteristic vector field of $\theta$) Consider the equation for
an expanding soliton:

\begin{eqnarray}
&f_{,0}=-W-\frac{1}{t}.&
\end{eqnarray}
 
Substituting the commutation relation $if_{,0}=f_{,1{\bar 1}}-
f_{,{\bar 1}1}$ in (3.3) and differentiating it in $Z_{\bar 1}$
and $Z_1$ directions, we get

\begin{eqnarray}
&f_{,1{\bar 1}1{\bar 1}}-f_{,{\bar 1}11{\bar 1}}=-iW_{,1{\bar 1}}.&
\end{eqnarray}

On the other hand, $X_f$ being a $CR$ vector field means ${\cal L}_
{X_f}J=0$, which is equivalent to (e.g.,[CL])

\begin{eqnarray}
&f_{,11}+iA_{11}f=0.&
\end{eqnarray}

Differentiating (3.5) in the $Z_{\bar 1}$ direction and exchanging
$1$ and $\bar 1$ using the commutation relation ([L2]), we obtain

\begin{eqnarray}
&f_{,1{\bar 1}1}+if_{,10}+f_{,1}W+iA_{11,{\bar 1}}f+iA_{11}f_{,{\bar 1}}
=0.&
\end{eqnarray}

Differentiating (3.3) in the $Z_1$ direction and switching $1$ and $0$
give $f_{,10}+A_{11}f_{,{\bar 1}}=-W_{,1}$. Substituting this in (3.6),
we obtain

\begin{eqnarray}
&f_{,1{\bar 1}1}-iW_{,1}+f_{,1}W+iA_{11,{\bar 1}}f=0.&
\end{eqnarray}

Also differentiating (3.3) twice in the $Z_{\bar 1}$ and $Z_1$
directions and exchanging ${\bar 1}1$ and $0$ using the commutation
relations, we obtain

\begin{eqnarray}
&f_{,{\bar 1}10}=-W_{,{\bar 1}1}-A_{11}f_{,{\bar 1}{\bar 1}}-
A_{11,{\bar 1}}f_{,{\bar 1}}-A_{{\bar 1}{\bar 1}}f_{,11}-
A_{{\bar 1}{\bar 1},1}f_{,1}.&
\end{eqnarray}

Differentiating the complex conjugate of (3.7) in the $Z_1$ direction
and substituting the result and (3.8) in the commutation relation:
$f_{,{\bar 1}11{\bar 1}}=f_{,{\bar 1}1{\bar 1}1}+if_{,{\bar 1}10}$,
we obtain an expression for $f_{,{\bar 1}11{\bar 1}}$. Substituting
this for the second term and the result of differentiating
(3.7) in the $Z_{\bar 1}$ direction for the first term in (3.4),
we can reduce (3.4) to

\begin{eqnarray}
&2i(W_{,1{\bar 1}}+W_{,{\bar 1}1})-(f_{,1{\bar 1}}-f_{,{\bar 1}1})W
+f_{,{\bar 1}}W_{,1}-f_{,1}W_{,{\bar 1}}&\nonumber\\
&+iA_{11}f_{,{\bar 1}{\bar 1}}+iA_{{\bar 1}{\bar 1}}f_{,11}
-i(A_{11,{\bar 1}{\bar 1}}+A_{{\bar 1}{\bar 1},11})f=0.&
\end{eqnarray}

Substituting $f_{,1{\bar 1}}-f_{,{\bar 1}1}=if_{,0}=-iW-it^{-1}$ (by
(3.3)) in (3.9), using (3.5) to replace $f_{,11}$ ($f_{,{\bar 1}{\bar 1}}$,
resp.) by $-iA_{11}f$ ($iA_{{\bar 1}{\bar 1}}f$, resp.), and noticing
the Bianchi identity: $A_{11,{\bar 1}{\bar 1}}+A_{{\bar 1}{\bar 1},11}
=W_{,0}$ and the definition of the sublaplacian operator $\Delta_b$, we finally obtain

\begin{eqnarray}
&2{\Delta_b}W+<{\nabla}W, X_{f}>+W^{2}+\frac{W}{t}=0&
\end{eqnarray}

\noindent in which ${\nabla}W={\nabla_b}W+W_{,0}T$ and $<{\nabla}W, X_{f}>
=-W_{,0}f-<{\nabla_b}W,J({\nabla_b}f)>_{J,{\theta}}$.

In the Riemannian case ([Ch]), we can add a certain quadratic term in
the involved vector field to get the Harnack quantity. However, in
our case, adding a quadratic term like $(constant)W|X_{f}|^{2}$ does not
seem to work without extra estimates on the torsion. So, as a first try,  
we assume the torsion vanishes at the initial time. It turns out that
the torsion vanishes for all time if our $CR$ structure $J$ is spherical.

\medskip

$\bf{Lemma\:3.1}$. Suppose $J$ is spherical and $A_{11}=0$ for an initial 
$\theta_{(0)}$. Then, under the Yamabe flow (1.1), $A_{11}$ vanishes for all $\theta_{(t)}$.

\medskip

$\bf{Proof}$: First, recall that the Cartan curvature tensor $Q_{11}$ is
related to the Tanaka-Webster curvature $W$ and torsion $A_{11}$ in the 
following formula: (Lemma 2.2 in [CL])

\begin{eqnarray}
Q_{11}=\frac{1}{6}W_{,11}+\frac{i}{2}WA_{11}-A_{11,0}-\frac{2i}{3}
A_{11,{\bar 1}1}.
\end{eqnarray}

The fundamental theorem of 3-dimensional $CR$ geometry due to Elie Cartan ([Ca]) asserts that $J$ being spherical is equivalent to $Q_{11}=0$. Now look at the evolution equation (2.8) of $A_{11}$. Each term in the right side of (2.8) contains the torsion or one of its derivatives in
view of (3.11). So obviously $A_{11}=0$ for all $t$ is a solution to (2.8). Therefore it suffices to show the uniqueness of solutions to
(2.8). However, using the commutation relation, we can write the right 
side of (2.8) as

\begin{eqnarray}
& 4(A_{11,1{\bar 1}}+A_{11,{\bar 1}1}-4iA_{11,0})-12WA_{11}. & \nonumber
\end{eqnarray}

\noindent The highest "weight" term is just $-4$ times the generalized Folland-Stein
operator ${\cal L}_{\alpha}$ (defined in [CL]) acting on $A_{11}$ with ${\alpha}=4$. Since ${\alpha}=4$ is not an odd integer, 
$-{\cal L}_{4}$ is subelliptic. So the uniqueness follows from the
standard theory for subparabolic equations. 

\begin{flushright}
Q.E.D.
\end{flushright}

\medskip

When the torsion $A_{11}$ vanishes identically, so does $W_{,0}$ due to
the Bianchi identity: ([L2]) 

\begin{eqnarray}
&A_{11,{\bar 1}{\bar 1}}+A_{{\bar 1}{\bar 1},11}=W_{,0}.&
\end{eqnarray}

\noindent Thus we can reduce $<{\nabla}W,X_{f}>$ in (3.10)  
to $-<{\nabla_b}W,J({\nabla_b}f)>_{J,{\theta}}$. Note that the vector
field ${\eta}=-J({\nabla_b}f)$ belongs to $\xi$, the contact bundle,
at each point. We call such a vector field a Legendrian vector field.
Releasing the $f$-dependence, we therefore consider (1.3) for arbitrary Legendrian vector field $\eta$ as our "Harnack" quantity.
(the coefficient $\frac{1}{8}$ of the last term in (1.3) is the minimal
value for (1.4) to hold as we'll see in the proof of Theorem A)
  
\bigskip

\section{\bf{The Harnack inequality: Proof of Theorems}}
\setcounter{equation}{0}

To apply the maximum principle to $Z({\theta},{\eta})$, we compute
the evolution equation for $Z({\theta},{\eta})$. For convenience,
we define ${\Box}={\partial}_{t}-4{\Delta_b}$ and ${\Box}{\eta}
=({\Box}{\eta_1})Z_{\bar 1}+({\Box}{\eta_{\bar 1}})Z_{1}$ for a
Legendrian vector field ${\eta}={\eta_1}Z_{\bar 1}+{\eta_{\bar 1}}Z_{1}$,
in which ${\Box}{\eta_1}={\partial_t}{\eta_1}-4({\eta}_{1,1{\bar 1}}+
{\eta}_{1,{\bar 1}1}$. Also we define the modulus of "Legendrian 2-tensor"
${\nabla_b}{\eta}$ as follows: $|{\nabla_b}{\eta}|^{2}_{J,{\theta}}=
2({\eta}_{1,{\bar 1}}{\eta}_{{\bar 1},1}+{\eta}_{1,1}{\eta}_{{\bar 1},
{\bar 1}})$. (recall that $h_{1{\bar 1}}=1$ and we express all tensors
using subindices)

\medskip
 
$\bf{Lemma\:4.1}$. Under the Yamabe flow (1.1), we have the following
evolution equations:

\begin{eqnarray}
&{\Box}(2{\Delta_b}W+W^{2})=12W{\Delta_b}W-4|{\nabla_b}W|_{J,{\theta}}^{2}+4W^{3}&\\
&{\Box}(W/t)=2W^{2}/t-W/t^{2}&\\
&{\Box}(W|{\eta}|^{2}_{J,{\theta}})=2W^{2}|{\eta}|^{2}_{J,{\theta}}
+2W<{\eta},{\Box}{\eta}>_{J,{\theta}}-8W|{\nabla_b}{\eta}|^{2}_{J,{\theta}}
-8<{\nabla_b}W,{\nabla_b}(|{\eta}|^{2}_{J,{\theta}})>_{J,{\theta}}&
\end{eqnarray}

\medskip

$\bf{Proof}$: (4.2) follows from (2.7). (4.1) follows from (2.12) with
$f=W$, (2.7), and the product formula: ${\Delta_b}(fg)=({\Delta_b}f)g+
f{\Delta_b}g+2<{\nabla_b}f,{\nabla_b}g>_{J,{\theta}}$. 
Similarly, (4.3) follows from a direct computation using (2.7) and the
above product formula. (noting that $|{\eta}|^{2}_{J,{\theta}}=2{\eta_1}
{\eta_{\bar 1}}$)

\begin{flushright}
Q.E.D.
\end{flushright}

\medskip

Let $\eta^1$, $\eta^2$
be two 2-tensors with components ${\eta}^{1}_{cd},{\eta}^{2}_{cd}$,respectively, in which $c,d=1$ or $\bar 1$. We define the Levi form for $\eta^1$, $\eta^2$
: $<{\eta^1},{\eta^2}>_{J,{\theta}}={\Sigma}{\eta}^{1}_{cd}{\eta}^{2}_
{{\bar c}{\bar d}}$, in which the sum is taken for all possible $(c,d)$
and ${\bar {({\bar 1})}}=1$ by convention. We define
${\nabla_b}{\eta}$ to be a 2-tensor with components ${\eta_{c,d}}$ 
for $\eta_c$ being components of a Legendrian vector field $\eta$.
Then the modulus
of ${\nabla_b}{\eta}$ defined previously is just the square root of the
above Levi form for ${\eta^1}={\eta^2}={\nabla_b}{\eta}$. 

\medskip
 
$\bf{Lemma\:4.2}$. Suppose the torsion $A_{11}$ vanishes identically
under the Yamabe flow (1.1). Then we have the following evolution
equation:

\begin{eqnarray}
&{\Box}(<{\nabla_b}W,{\eta}>_{J,{\theta}})=W<{\nabla_b}W,{\eta}>_{J,
{\theta}}+<{\nabla_b}W,{\Box}{\eta}>_{J,{\theta}}-8<{\nabla_b}({\nabla_b}W),{\nabla_b}{\eta}>_{J,{\theta}}&
\end{eqnarray}

\medskip

$\bf{Proof}$: First observe that 

\begin{eqnarray}
{\Box}(W_{,1})&=&{\partial_t}(Z_{1}W)-4(W_{,11{\bar 1}}+W_{,1{\bar 1}1})\\
              &=&5WW_{,1}+4(W_{,{\bar 1}11}-W_{,11{\bar 1}})(by\:(2.7)\:
                             and\:{\dot Z}_{1}=WZ_{1})\nonumber\\
&=&WW_{,1}-8iW_{,01}\:mod(A_{11},A_{11,{\bar 1}})\:(by\:commutation\:
   relations)\nonumber\\
&=&WW_{,1}.\:(A_{11}=0\:and\:W_{,0}=0\:by\:(3.12))\nonumber
\end{eqnarray}

Now (4.4) follows from a direct computation using the product formula and (4.5).

\begin{flushright}
Q.E.D.
\end{flushright}

\medskip

By combining Lemmas 4.1 and 4.2, we find that if the torsion vanishes
identically, the evolution equation
for $Z({\theta},{\eta})$ is given by

\begin{eqnarray}
{\Box}Z({\theta},{\eta})&=&12W{\Delta_b}W-4|{\nabla_b}W|_{J,{\theta}}^{2}+4W^{3}+2W^{2}/t-W/t^{2}\\
&&+W<{\nabla_b}W,{\eta}>_{J,{\theta}}+\frac{1}{4}W^{2}|{\eta}|^{2}_{J,{\theta}}+<{\nabla_b}W+
\frac{1}{4}W{\eta},{\Box}{\eta}>_{J,{\theta}}\nonumber\\
&&-8<{\nabla_b}({\nabla_b}W),{\nabla_b}{\eta}>_{J,{\theta}}
-W|{\nabla_b}{\eta}|^{2}_{J,{\theta}}
-<{\nabla_b}W,{\nabla_b}(|{\eta}|^{2}_{J,{\theta}})>_{J,{\theta}}.\nonumber
\end{eqnarray}

$\bf{Proof\:of\:Theorem\:A}$: First observe that by Corollary 2.2, $W$ is 
always positive, and for all $\eta$, $Z({\theta},{\eta}){\geq}
Y({\theta})$ where

\begin{eqnarray}
&Y({\theta})=2{\Delta_b}W+W^{2}+\frac{W}{t}-2W^{-1}|{\nabla_b}W|^
{2}_{J,{\theta}}.&
\end{eqnarray}

\noindent Also there exists a positive constant $\delta$ such that $
Y({\theta})>0$ for $t<{\delta}$, hence $Z({\theta},{\eta})>0$ for
$t<{\delta}$.

Suppose $Z({\theta},{\eta}){\leq}0$ at some space-time point for some
$\eta$. Then there exists a first time ${\tau}>0$, a point ${\zeta}{\in}
M$ and a Legendrian tangent vector $\eta$ at $\zeta$ such that
at $({\zeta},{\tau})$,

\begin{eqnarray}
&Z({\theta},{\eta})=0.& 
\end{eqnarray}

We extend $\eta$ so that at $({\zeta},{\tau})$, 

\begin{eqnarray}
&{\eta_{1,{\bar 1}}}=W^{-1}(4iW_{,1{\bar 1}}-W_{,{\bar 1}}{\eta_1})&\\
&{\eta_{1,1}}=-W^{-1}W_{,1}{\eta_1}&
\end{eqnarray}

\noindent where ${\eta}=i{\eta_1}Z_{\bar 1}-i{\eta_{\bar 1}}Z_{1}$.
Substituting (4.9), (4.10) in the last three terms of (4.6), involving first derivatives of $\eta$, we obtain that at $({\zeta},{\tau})$, 

\begin{eqnarray}
&-8<{\nabla_b}({\nabla_b}W),{\nabla_b}{\eta}>_{J,{\theta}}
-W|{\nabla_b}{\eta}|^{2}_{J,{\theta}}
-<{\nabla_b}W,{\nabla_b}(|{\eta}|^{2}_{J,{\theta}})>_{J,{\theta}}&
\nonumber\\ 
&=W^{-1}(2|4iW_{,1{\bar 1}}-W_{,{\bar 1}}{\eta_1}|^{2}+
2|W_{,{\bar 1}}{\eta_{\bar 1}}|^{2}).&
\end{eqnarray}

\noindent In deriving (4.11), we have used $W_{,11}=0$ for all time 
due to Lemma 3.1 and $Q_{11}=0$ in (3.11). (note that the choice of
${\nabla_b}{\eta}$ in (4.9),(4.10) is to maximize the left side of (4.11))

Now if ${\nabla_b}W+\frac{1}{4}W{\eta}{\neq}0$ at $({\zeta},{\tau})$,
we extend $\eta$ by choosing the value of ${\Box}{\eta}$ at $({\zeta},{\tau})$ to kill all terms on the right side of (4.6) except, say, the term
$2W^{2}/t$. Then it follows that at 
$({\zeta},{\tau})$,
$0{\geq}{\partial_t}Z=4{\Delta_b}Z+2W^{2}/{\tau}{\geq}2W^{2}/{\tau}$,
a contradiction. So we assume 

\begin{eqnarray}
&{\nabla_b}W+\frac{1}{4}W{\eta}=0&
\end{eqnarray}

\noindent at $({\zeta},{\tau})$. By (4.8) and (4.12), we can express
${\Delta_b}W$ in terms of $W$ and $\eta$ at $({\zeta},{\tau})$:

\begin{eqnarray}
&2{\Delta_b}W=\frac{1}{8}W|{\eta}|_{J,{\theta}}^{2}-W^{2}-\frac{W}{t}.&
\end{eqnarray}

From Lemma 3.1 and (3.12), $W_{,0}$ vanishes identically for all time.
It follows that $W_{,1{\bar 1}}=W_{,{\bar 1}1}$ by the commutation
relation. Therefore ${\Delta_b}W=2W_{,1{\bar 1}}$, so by (4.13) we
can express $W_{,1{\bar 1}}$ at $({\zeta},{\tau})$ as follows:

\begin{eqnarray}
&W_{,1{\bar 1}}=\frac{1}{4}(\frac{1}{8}W|{\eta}|^{2}_{J,{\theta}}-
W^{2}-\frac{W}{t}).&
\end{eqnarray}

Now substituting (4.14),(4.12) in (4.11) and making use of (4.13),(4.12),
we can reduce (4.6) to an expression in $W$,$\eta$ only:

\begin{eqnarray}
&{\Box}Z=\frac{W}{t^{2}}+\frac{1}{2}W^{2}|{\eta}|^{2}_   
{J,{\theta}}+\frac{1}{32}W|{\eta}|_{J,{\theta}}^{4}&\nonumber
\end{eqnarray}

\noindent at $({\zeta},{\tau})$. Hence the maximum principle implies that
at $({\zeta},{\tau})$,

\begin{eqnarray}
&0{\geq}{\partial_t}Z=4{\Delta_b}Z+\frac{W}{t^{2}}+\frac{1}{2}W^{2}|{\eta}|^{2}_   
{J,{\theta}}&\nonumber\\
&+\frac{1}{32}W|{\eta}|_{J,{\theta}}^{4}{\geq}W{\tau}^{-2},&\nonumber
\end{eqnarray}

\noindent which is a contradiction (to Corollary 2.2).
So $Z({\theta},{\eta})>0$ completing the proof of Theorem A.

\begin{flushright}
Q.E.D.
\end{flushright}

Note that we actually obtain the strict inequality from the proof of
Theorem A. Taking ${\eta}=-4W^{-1}{\nabla_b}W$ in Theorem A implies that

\begin{eqnarray}
&Z({\theta},{\eta})=Y({\theta}){\geq}0,& 
\end{eqnarray}

\noindent where $Y({\theta})$ is defined in (4.7). 

$\bf{Proof\:of\:Theorem\:B}$: By (2.7) we rewrite (4.15) as

\begin{eqnarray}
&{\partial_t}W+\frac{2W}{t}-4W^{-1}|{\nabla_b}W|^{2}_{J,{\theta}}{\geq}0.&
\end{eqnarray}

Integrating (4.16) over Legendrian paths connecting $x_1$,$x_2$ from
$t_1$ to $t_2$ and using $az+{\bar a}{\bar z}+b|z|^{2}{\geq}-b^{-1}
|a|^{2}$ for $b>0$, we obtain (1.5).

\begin{flushright}
Q.E.D.
\end{flushright}

\bigskip

\section{\bf{Nontriviality of initial conditions}}
\setcounter{equation}{0}

In this section, we will construct contact forms on the standard $CR$
3-sphere $(S^{3},{\hat {\xi}},{\hat J})$ with nonconstant positive
Tanaka-Webster curvature and vanishing torsion.

Suppose our $S^{3}$ is
defined by $|z_{1}|^{2}+|z_{2}|^{2}=1$ for $(z_{1},z_{2}){\in}C^{2}$.
The standard contact form ${\hat{\theta}}=i({\sigma}-{\bar{\sigma}})$
where ${\sigma}=z_{1}d{\bar z}_{1}+z_{2}d{\bar z}_{2}$. Take ${\hat
{\theta}}^{1}=\sqrt{2}(z_{1}dz_{2}-z_{2}dz_{1})$ such that (2.1) is
satisfied for the "hatted" quantities. It is easy to deduce ${\hat Z}_{1}
=\frac{1}{\sqrt{2}}({\bar z}_{1}{\partial_2}-{\bar z}_{2}{\partial_1})$
where ${\partial_j}=\frac{\partial}{{\partial}z_{j}}$ for $j=1,2$.
Also ${{\hat{\omega}}_{1}}^{1}=-2({\bar z}_{1}dz_{1}+{\bar z}_{2}dz_{2})$
and ${\hat A}_{11}=0$ by (2.2).

Now  let ${\theta}=e^{2{\lambda}}{\hat{\theta}}$. It follows from (2.5)
that $A_{11}=0$ (with respect to $\theta$) if and only if 

\begin{eqnarray}
&{\lambda}_{,11}=2({\lambda}_{,1})^{2}&
\end{eqnarray}

\noindent in which ${\lambda}_{,1}={\hat Z}_{1}{\lambda}$ and ${\lambda}_{,11}=({\hat Z}_{1})^{2}{\lambda}$ since ${{\hat{\omega}}_{1}}^{1}({\hat Z}_{1})=0$. It is a direct verification that ${\lambda}=-ln|az_{1}+bz_{2}+c|$ (well defined on $S^3$ for $|c|{\gg}|a|,|b|$) satisfies (5.1).
Next we compute $W$ from the formula (2.6) for ${\hat W}=1$ and $\lambda$
given above. The final result is

\begin{eqnarray}
&W=-(3|z_{1}|^{2}+2|z_{2}|^{2})|a|^{2}-(3|z_{2}|^{2}+2|z_{1}|^{2})|b|^{2}
-(a{\bar b}z_{1}{\bar z}_{2}+{\bar a}b{\bar z}_{1}z_{2})&\\
&-c({\bar a}{\bar z}_{1}+{\bar b}{\bar z}_{2})-{\bar c}(az_{1}+bz_{2})
+|c|^{2}.&\nonumber
\end{eqnarray}

Now it is easy to see from (5.2) that $W$ is positive on $S^3$ for $|c|{\gg}|a|,|b|$,
and nonconstant in general.

\bigskip

\begin{tabular}{ll}
Chang: Department of Mathematics & Cheng: Institute of Mathematics \\
National Tsing-Hua University & Academia Sinica, Nankang \\
Hsinchu, Taiwan, R.O.C. & Taipei, Taiwan, R.O.C. \\
E-mail: scchang{\bf @}math.nthu.edu.tw & E-mail: cheng{\bf@}math.sinica.edu.tw
\end{tabular}



\bigskip

\begin{thebibliography}{W}


\bibitem[Ca]{ca}
E. Cartan,
{\em Sur la g\'{e}ometrie pseudo-conforme des hypersurfaces de 
l'espace de deux variables complexe},
I, Ann. Mat., 11(1932),17-90; II, Ann. Sc. Norm. Sup. Pisa, 1(1932),333-354


\bibitem[Ch]{ch}
B. Chow,
{\em The Yamabe flow on locally conformally flat manifolds with
positive Ricci curvature},
Commun. Pure and Appl. Math., XLV(1992),1003-1014

\bibitem[CL]{cl}
J.-H. Cheng and J. M. Lee,
{\em The Burns-Epstein invariant and deformation of CR structures},
Duke Math. J., 60(1990),221-254

\bibitem[H1]{h1}
R. S. Hamilton,
{\em The Harnack estimate for the Ricci flow},
J. Diff. Geom., 37(1993),225-243 

\bibitem[H2]{h2}
--------,
{\em The Ricci flow on surfaces},
Math. and General Relativity, Contemporary Math. 71,(1988),237-262


\bibitem[JL1]{jl1}
D. Jerison and J. M. Lee,
{\em The Yamabe problem on $CR$ manifolds},
J. Diff. Geom., 25(1987),167-197

\bibitem[JL2]{jl2}
-----------,
{\em Extremals for the Sobolev inequality on the Heisenberg group
and the $CR$ Yamabe problem},
J. Amer. Math. Soc., 1(1988),1-13

\bibitem[JL3]{jl3}
-----------,
{\em Intrinsic $CR$ normal coordinates and the $CR$ Yamabe problem},
J. Diff. Geom., 29(1989),303-343

\bibitem[L1]{l1}
J. M. Lee,
{\em The Fefferman metric and pseudohermitian invariants},
Trans. Amer. Math. Soc., 296(1986),411-429

\bibitem[L2]{l2}
---------,
{\em Pseudo-Einstein structures on $CR$ manifolds},
Am. J. Math., 110(1988),157-178

\bibitem[Ta]{ta}
N. Tanaka,
{\em A Differential Geometric Study on Strongly
Pseudo-Convex Manifolds},
1975, Kinokuniya Co. Ltd., Tokyo

\bibitem[We]{we}
S. M. Webster,
{\em Pseudohermitian structures on a real hypersurface},
J. Diff. Geom., 13(1978),25-41

\end{thebibliography}
\end{document}